\documentclass[12pt]{article}

\usepackage{amsmath, amsfonts, amssymb, latexsym}

\usepackage{comment}
\usepackage{rotating}

\usepackage{color}

\usepackage{amsthm}

\newtheorem{theorem}{Theorem}

\theoremstyle{definition}

\theoremstyle{plain}
\theoremstyle{remark}
\newtheorem{remark}{\textbf{Remark}}

\begin{document}
\title{On linear forms containing the Euler constant\thanks{ The work was supported by Programm 1 DMS RAS and grants N-Sh-3906.2008.1,
 RFBR-08-01-00179.}
}
\author{A.I.Aptekarev}

\date{}
\maketitle

\setcounter{equation}{0}
\section{Introduction}\label{S:Intro}
We study the arithmetical nature of two numbers
$$
\gamma:=-\int_0^\infty \ln x \,e^{-x} dx \qquad \mbox{and} \qquad \delta:=\int_0^\infty \ln (x+1) \,e^{-x} dx\,.
$$
The first number is the famous Euler (or Euler-Mascheroni) constant
$$
\gamma:=\lim\limits_{n\to\infty}
\left(\sum\limits_{k=1}^n\frac{1}{k}-\ln(n)\right)\,\approx\,0.577\ldots\;,
$$
The second number is called Euler-Gompertz constant (see \cite{FI}).
A relation
$$
\delta:=\int_0^\infty \frac{e^{-x} dx}{x+1}\,\approx\,0.596\ldots
$$
to the Laguerre polynomials gives a sequence of rational  approximants for~$\delta$
\begin{equation}\label{1}
\frac{\tilde{p}_n}{\tilde{q}_n}\rightarrow\delta\;,\;n\to\infty\;,
\end{equation}
generated by the recurrence relations
$$
\tilde{q}_{n+1}\,\,=\,2\,(n+1)\,\tilde{q}_n\,-\,n^2\,\tilde{q}_{n-1}\,,
$$
with initial condition
\begin{equation*}
\begin{array}{l}
\tilde{p}_0=0\,,\quad\tilde{p}_1=1\,,\\
\tilde{q}_0=1\,,\quad\tilde{q}_1=2\,.\\
\end{array}
\end{equation*}
The Perron asymptotics for the Laguerre polynomials
\begin{equation*}
\begin{array}{l}
\tilde{q}_n\,=\,n!\,\displaystyle\frac{e^{2\sqrt{n}}}{\sqrt[4]{n}}\,\left(\frac{1}{2\sqrt{\pi e}}+
O(n^{-1/2})\right)\;,\\
\\
\tilde{p}_n-\delta \tilde{q}_n\,=\,O\left(n!\,\displaystyle\frac{e^{-2\sqrt{n}}}{\sqrt[4]{n}}\right)\;,\\
\end{array}
\end{equation*}
confirm \eqref{1}. However,  the rational approximants \eqref{1} do not imply that
$\delta$ is an irrational number. As far as we know, irrationality  of $\delta$ is
still an open problem, as well as the famous open problem about irrationality of
the Euler constant $\gamma$.

\section{Result}\label{S:Res}

\begin{theorem}\label{T2}
Given  sequences of numbers
$$
u:=\{u_n\}\,,\quad v:=\{v_n\}\,,\quad w:=\{w_n\}\,,
$$
generated by the recurrence relations
\begin{equation}\label{2}
\begin{array}{l}
(16n+1)(16n-15)\,u_{n+1}\,\,=\,(16n-15)(256n^3+528n^2+352n+73)\,u_n
\\\phantom{--uuuuuuuuuuuuuuuuuuu-}-\,(16n+17)(128n^3+40n^2-82n-45)\,u_{n-1}
\\\phantom{--uuuuuuuuuuuuuuuuuuuuuuuuuuuuuuuuuu-}+\,n^2(16n+17)(16n+1)\,u_{n-2}\;,
\end{array}
\end{equation}
with  initial conditions
\begin{equation*}
\begin{array}{l}
u_0=-2\,,\;\!u_1=\phantom{-2}\!7\,,\;u_2=\phantom{-1}558\,,\\
v_0=-1\,,\;v_1=-22\,,\;v_2=-1518\,,\\
w_0=\phantom{-}\!\!0\,,\;w_1=-17\,,\;\!w_2=-1209\,,\\
\end{array}
\end{equation*}
Then

1) $\quad u_n,\,v_n\,, w_n\,\in \mathbb{Z}$;

2) difference equation \eqref{2} has  solutions with three different asymptotics
as $n \rightarrow \infty$

\begin{equation}\label{13}
\begin{array}{l}
u_n\,,\,\,v_n\,,\,\,w_n\,\,=\,\,O\left(\displaystyle\frac{(2n)!\,4^n}{n^{3/2}}\right)\;,\\
\\
\left[w_n +(e \gamma +\delta) u_n\right],\,\,\left[v_n +e u_n\right]\,\,=\,\,O\left(\displaystyle\frac{e^{\sqrt{2n}}\,n^{5/4}}{4^n}\right)\;,\\
\\
l_n:=u_n \delta-\gamma v_n + w_n\,=\,\displaystyle\frac{n^{5/4}}{e^{\sqrt{2n}}\,4^n}\,
\left(\frac{2\sqrt[4]{2}}{e^{3/8}}+O(n^{-1/2})\right)\;.\\
\end{array}
\end{equation}
 \end{theorem}

\bigskip

\begin{remark}
The asymptotics of $l_n$ in \eqref{13} give a quantitative characterization of the
fact that one of two constants $\gamma$ or $\delta$ is an irrational number. The
validity of this
fact was known before, it follows from the A.B.~Shidlovski result~\cite{Shi} on
algebraic
independence of the values of $E$-functions (see also K.~Mahler \cite{Mah}).
Indeed
(see, for example, the last statement of~\cite{Mah}),  numbers
$$
1\,-\,\frac{1}{e}\quad \mbox{and} \quad -\left(\gamma\,-\,\frac{\delta}{e}\right)
$$
are algebraically independent. It implies that numbers $\gamma$ and $\delta$
can not be  rational numbers simultaneously.
\end{remark}

\medskip

We
thank Yu.V.~Nesterenko and T.~Rivoal for indication  of this implication to us.

\section{A rational approximation of Euler constant}\label{S:RAG}
The numbers generated by \eqref{2} are intimately related with  rational
approximants of the Euler constant (studied in \cite{A1}). These approximants are
produced by the functional Hermite-Pade rational approximants for a system of four
functions $\{\hat{\mu}_1, \, \hat{\mu}_2,\, \hat{s}_1, \, \hat{s}_2 \}$~:
\begin{equation}\label{fun}
\hat{\mu}_k (z) \,: = \, \int_0^1\frac{w_k(z) dx}{z - x} \, , \qquad \hat{s}_k (z)
\,: = \, \int_1^\infty \frac{w_k(z) dx}{z - x} \, ,
\end{equation}
where
\begin{equation}\label{weigh}
w_k(x) \,: = \, x^{\alpha_k} \, (1-x)^\alpha \, e^{-\beta x} \, , \qquad k=1,2 \,.
\end{equation}
For the denominators   $Q_n$ of the Hermite-Pade rational approximants of this
system \eqref{fun}-\eqref{weigh} the generalized Rodrigues formula was known (see
the example of subsection 2.2 in~\cite{ABVA}; for similar systems, see
also~\cite{S1} and ~\cite{S2})
\begin{equation}
Q_{n}(z)\,=\,\frac{1}{(n!)^2} \, \, w_2^{-1}{d^n \over dz^n}\,
 \left[w_2 z^n w_1^{-1} {d^n \over dz^n} \left[w_1 z^n (1-z)^{2n}\right] \right]
 \,\,.
\label{rodr}
\end{equation}
Then the above mentioned rational approximants $\displaystyle \frac{p_n}{q_n}$ of
the Euler constant $\gamma$ are defined as
\begin{equation}
p_n \, - \, \gamma\, q_n \,:= \, \int_0^\infty Q_n(x) \ln (x)\, e^{-x}\, dx \,=:\,
f_n , \label{form}
\end{equation}
where  $Q_n$ is taken from~\eqref{rodr} -- \eqref{weigh} with parameters
$\alpha_1=\alpha_2=0$ ¨ $\alpha = - \beta =1$.

\medskip

In \cite{X}--\cite{B} recurrence relations (of six, seven and eight terms) for
polynomials~\eqref{rodr} were studied, which in~\cite{AT1} eventually brought a
four-term recurrence relation for sequences of numbers $p := \{ p_n \}$ and $q :=
\{ q_n \}$
\begin{equation}\label{8}
\begin{array}{l}
(16n-15)\,q_{n+1}\,\,=\,(128n^3+40n^2-82n-45)\,q_n\phantom{----uuuuuuuuu-}
\\\phantom{uuuu-}-\,n^2(256n^3-240n^2+64n-7)\,q_{n-1}\,+\,n^2(n-1)^2(16n+1)
\,q_{n-2}\;,
\end{array}
\end{equation}
with initial conditions
\begin{equation*}
\begin{array}{l}
p_0=0\,,\;p_1=2\,,\;p_2=31\,,\\
q_0=1\,,\;q_1=3\,,\;q_2=50\,.\\
\end{array}
\end{equation*}
We also highlight a solution $r := \{ r_n \}$ of difference equation~\eqref{8}
with initial conditions
\begin{equation*}
\begin{array}{l}
r_0=0\,,\;r_1=1\,,\;r_2=24\,.\\
\end{array}
\end{equation*}
The fact that  numbers $\,\,q_n,\,p_n,\,r_n\,\in \mathbb{Z}$ are integers for
$\,n\,\in \mathbb{N}$ was proven in~\cite{T2}. Finally in~\cite{AL} and~\cite{T1},
the following asymptotics for $\,\{ q_n,\,p_n,\,r_n\}\,$ and the linear forms with
these coefficients were obtained
\begin{equation}\label{9}
\begin{array}{rl}
q_n\,,\,\,p_n\,,\,\,r_n\,\,&=\,\,O\,(\displaystyle\frac{(2n)!\,e^{\sqrt{2n}}}{\sqrt[4]{n}})\;,\\
\\
f_n:\,=\,p_n-\gamma q_n\,&=
\,(2n)!\,\displaystyle\frac{e^{-\sqrt{2n}}}{\sqrt[4]{n}}\,
\left(\frac{2\sqrt{\pi}}{(4e)^{3/8}}+O(n^{-1/2})\right)\;,\\
\\
g_n:=e p_n-( e \gamma +\delta) q_n +r_n\,&=\,\,\displaystyle\frac{1}{16^{n}}\,
\left(\frac{1}{8}+O(n^{-1})\right)\;.\\
\end{array}
\end{equation}

\section{Proof of Theorem~\ref{T2}}\label{S:Proof}
\textbf{A)} We define
\begin{equation}\label{10}
u_n\,:=\, \displaystyle \frac{\Delta_n^{(qp)}}{(n!)^2}\,, \qquad v_n\,:=\,
\displaystyle \frac{\Delta_n^{(qr)}}{(n!)^2}\,,\qquad w_n\,:=\, \displaystyle
\frac{\Delta_n^{(p r)}}{(n!)^2}\,,
\end{equation}
where we use a notation
\begin{equation}\label{11}
\Delta_n^{(a b)}\,:=\, \det\left(\!\begin{array}{c}\!\!a_{n+1}\,,\;a_{n}\!\!
\\\!\!b_{n+1}\,,\;b_{n}\!\!\\\end{array}\right)\,,
\qquad a := \{ a_n \}\,,\quad b := \{ b_n \}\,.
\end{equation}
\textbf{B)} Substituting the recurrence relations
$$
\left(\!\begin{array}{c}\!\,a_{n+1}\!\!\\\!\,b_{n+1}\!\!\\\end{array}\right)\,= \,
A_n\,\!\left(\!\begin{array}{c}\!\,a_{n}\!\!\\\!\,b_{n}\!\!\\\end{array}\right)\,
+\,B_n\,\!\left(\!\begin{array}{c}\!\,a_{n-1}\!\!\\\!\,b_{n-1}\!\!\\\end{array}\right)\,
+\,C_n\,\!\left(\!\begin{array}{c}\!\,a_{n-2}\!\!\\\!\,b_{n-2}\!\!\\\end{array}\right)\,
$$
into  determinant~\eqref{11}, we get
$$
\left|\!\begin{array}{c}\!\,a_{n+1}\,,\;a_{n}\!\!\\\!\,b_{n+1}\,,\;b_{n}\!\!\\\end{array}\right|\,=
\,
-B_n\,\!\left|\!\begin{array}{c}\!\,a_{n}\,,\;a_{n-1}\!\!\\\!\,b_{n}\,,\;b_{n-1}\!\!\\\end{array}\right|\,
-\,C_n
A_{n-1}\,\!\left|\!\begin{array}{c}\!\,a_{n-1}\,,\;a_{n-2}\!\!\\\!\,b_{n-1}\,,\;b_{n-2}\!\!\\\end{array}\right|\,
+\,C_nC_{n-1}\,\!\left|\!\begin{array}{c}\!\,a_{n-2}\,,\;a_{n-3}\!\!\\\!\,b_{n-2}\,,\;b_{n-3}\!\!\\\end{array}\right|\,
\,,$$ that from~\eqref{8} leads to~\eqref{2}.

$\phantom{a}$\\
\textbf{C)} The fact that $\displaystyle
\frac{\Delta_n^{(qp)}}{(n!)^2}\,\in\,\mathbb{Z}$ is proven in~\cite{AT2}. We get
that $v_n,\,w_n\,\in\,\mathbb{Z}$ analogously.

$\phantom{a}$\\
\textbf{D)} The asymptotics of
$$
\Delta_n^{(q p)}\,=\, \left|\!\begin{array}{l}\!\,q_{n+1} \phantom{-\gamma
q_{n+1}}\,\,\,\,,\;q_{n}\!\!
\\\!\,p_{n+1}-\gamma q_{n+1}\,,\;p_{n}-\gamma q_{n}\!\!\\\end{array}\right|\,
=\, O\left( \displaystyle \frac{(2n)!^2}{n^{3/2}}\right)
$$
was computed in~\cite{AT2}. The same way, from \eqref{9} we deduce
$$
\Delta_n^{(q r)}\,=\, \left|\!\begin{array}{l}\!\,q_{n+1} \phantom{-\gamma
q_{n+1}}\,\,\,\,,\;q_{n}\!\!
\\\!\,r_{n+1}-\gamma q_{n+1}\,,\;r_{n}-\gamma q_{n}\!\!\\\end{array}\right|\,
=\, O\left( \displaystyle \frac{(2n)!^2}{n^{3/2}}\right)\,.
$$
Noticing that
$$
\Delta_n^{(p r)} \frac{1}{e}\,-\,\Delta_n^{(p q)} \left( \gamma + \frac{\delta
}{e} \right)\,=\,\Delta_n^{(p g)}\,,
$$
we get from \eqref{9} asymptotics for $[ w_n + ( e \gamma + \delta)\, u_n ]$ and
$w_n$. Analogously, the identity
$$
\Delta_n^{(q p)} \,+\,\Delta_n^{(q r)} \frac{1 }{e} \,=\,\Delta_n^{(q g)}
$$
brings the asymptotics for $( v_n +  e \, u_n )$. Finally,
$$
- \Delta_n^{(p q)} \frac{\delta }{e} \,+\,
 \Delta_n^{(p r)} \frac{1}{e}\,-\,\Delta_n^{(q r)}  \frac{\gamma
}{e} \,=\,\Delta_n^{(f g)}
$$
from \eqref{9} brings the asymptotics for $l_n$ in \eqref{13}.

$\phantom{a}$\\
The Theorem is proved.

\vskip 4. cm

\begin{abstract}

We present linear forms with integer coefficients containing the Euler-Mascheroni
and Euler-Gompertz constants. The forms are defined by four-terms recurrence
relations. Asymptotics of the forms and their coefficients are obtained.

\end{abstract}

\end{document}